\documentclass{article} 
\usepackage{amssymb, amsfonts, amsmath, amsthm}
\usepackage{epsfig,lpic,wrapfig}
\usepackage[T2A,T1]{fontenc}
\usepackage[utf8]{inputenc}
\usepackage[russian,english]{babel}
\usepackage{graphicx}
\usepackage{url}
\usepackage[top=0.9in, bottom=0.9in,left=0.9in, right=0.9in, paperwidth=6in, paperheight=9in]{geometry}
\usepackage[colorlinks=true,
citecolor=black,
linkcolor=black,
anchorcolor=black,
filecolor=black,
menucolor=black,
urlcolor=black%
]{hyperref}
\hypersetup{pdftitle={Плоское оригами и построения.},
pdfauthor={Антон Петрунин}}
\usepackage{multicol}

\def\noi{\noindent}%

\def\phi{\varphi}%
\def\"#1{{\accent"7F #1\penalty10000\hskip 0pt plus 0pt}} 

\begin{document}
\otherlanguage{russian}
\hyphenation{Фу-дзи-та}
\title{Плоское оригами и построения}
\author{А.~Петрунин}
\date{}
\maketitle

В этой заметке мы обсудим геометрические построения связанные с оригами, то есть складыванием бумаги. 
Этот тип построений основан на таком наблюдении: 
если разгладить на столе сложенный лист бумаги то складка образует прямую линию.
В таких построениях можно обходится без карандаша, про прямую можно думать как про линию сгиба, а про точку как пересечение двух прямых. 

Например, если нам нужно построить срединный перепендикуляе к отрезку $[AB]$ то можно сложить лист так, чтобы эти точки $A$ и $B$ совпали и потом разгладить лист на столе.
Линия сгиба окажется искомым перпендикуляром. 

\begin{wrapfigure}{r}{20mm}
\begin{lpic}[t(-8mm),b(-2mm),r(0mm),l(0mm)]{pics/perpendicular-bisector(1)}
\lbl[r]{3,11.5;$A$}
\lbl[l]{16,11.5;$B$}
\end{lpic}
\end{wrapfigure}

Заметим, что нам ничего не потребовалось кроме стола и бумаги,
тогда как
классическое построение требует чертёжных инструментов и кроме того длиннее, см. рисунок.

Чтобы продолжить нам нужно описать првила оригами-построений.

Например в построениях с помощью циркуля и линейки 
разрешается отмечать на плоскости произвольную точку отличную от уже отмеченных и делать такие операции:
\begin{enumerate}
\item Провести прямую через две отмеченные точки;
\item Построить окружность с центром в отмеченной точке и радиусом равным расстоянию между отмеченными точками;
\item Отметить точки пересечения 
двух прямых 
или двух окружностей
или прямой с окружностью. 
\end{enumerate}

Ниже мы приведём 7 правил которыми пользуются оригамисты.
Подобно построению с помощью циркуля и линейки, 
оригами-построение состоит из последовательного складывания листа по правилам из этого списка.
Некоторые из этих складок  
можно получить как результат последовательного применения остальных.
Математик обошёлся бы без них, 
но оригамисты не мнут бумагу зря.

\newpage
\begin{multicols}{2}
\begin{enumerate}
\item\label{psps} 
Лист можно сложить так, что две отмеченные точки будут на складке.
\begin{center}
\begin{lpic}[t(-1mm),b(-1mm),r(0mm),l(0mm)]{pics/psps1(1)}
\end{lpic}
\end{center}

\item\label{ptop} Лист можно сложить так, что одна отмеченная точка перейдёт в другую отмеченную точку.
\begin{center}
\begin{lpic}[t(-1mm),b(-1mm),r(0mm),l(0mm)]{pics/ptop1(1)}
\end{lpic}
\end{center}
\item\label{ltol} Лист можно сложить так, что отмеченная прямая перейдёт в другую отмеченную прямую.
\begin{center}
\begin{lpic}[t(-1mm),b(-1mm),r(0mm),l(0mm)]{pics/ltol1(1)}
\end{lpic}
\end{center}
\item\label{psls} Лист можно сложить так, что отмеченная  точка попадёт на складку, а отмеченная прямая перейдёт в себя (то есть, линия складки будет ей перпендикулярна).
\begin{center}
\begin{lpic}[t(-1mm),b(-1mm),r(0mm),l(0mm)]{pics/psls1(1)}
\end{lpic}
\end{center}
\item\label{ptolps} Пусть отмеченны прямая $p$ и две точки $A$ и $B$.
Тогда лист можно\footnote{Такая складка существует не всегда, правило утверждает только, что если такая складка есть, то её «можно» найти.} сложить так, что точка $A$ попадёт на складку, а $B$ на прямую $p$.
\begin{center}
\begin{lpic}[t(-1mm),b(-1mm),r(0mm),l(0mm)]{pics/ptolps1(1)}
\lbl[t]{3,2;$A$}
\lbl[t]{14,2;$B$}
\lbl[b]{8,16;$p$}
\end{lpic}
\end{center}
\item\label{2ptol} Пусть отмечены две прямые $p$ и $q$ и две точки $A$ и $B$. 
Тогда лист можно$^2$ сложить так, что точка $A$ попадёт на прямую $p$, а точка $B$ попадёт на прямую $q$.
\begin{center}
\begin{lpic}[t(-1mm),b(-1mm),r(0mm),l(0mm)]{pics/2ptol1(1)}
\lbl[r]{.5,7;$A$}
\lbl[t]{6,1.5;$B$}
\lbl[l]{17,6;$q$}
\lbl[b]{5,18;$p$}
\end{lpic}
\end{center}
\item\label{ptolls} Пусть отмечены  две прямые $p$ и $q$ и точка $A$. 
Тогда лист можно сложить так, что точка $A$ попадёт на прямую $p$, а прямая $q$ перейдёт в себя (то есть, линия складки будет ей перпендикулярна).
\begin{center}
\begin{lpic}[t(-1mm),b(-1mm),r(0mm),l(0mm)]{pics/ptolls1(1)}
\lbl[r]{.5,7;$A$}
\lbl[b]{5,18;$p$}
\lbl[tl]{17,9;$q$}
\end{lpic}
\end{center}
\end{enumerate}
\end{multicols}
\newpage

Теперь можно заниматься построениями. 
Например:
\begin{enumerate}
\item Eсли задан треугольник,  его \textit{биссектрисы}, а, стало быть, и \textit{центр его вписанной окружности} можно найти, применив правило \ref{ltol} ко всем парам его сторон.

\item \textit{Срединные перпендикуляры} и \textit{центр описанной окружности} можно найти, применив правило \ref{ptop} ко всем парам его вершин.
После этого найти \textit{медианы} и \textit{центр тяжести}, применив правило \ref{psps} к каждой вершине в паре с уже найденной выше серединой противоположной стороны. 
\item \textit{Высоты} и \textit{ортоцентр} легче всего найти, применив правило \ref{psls} к каждой вершине в паре с противоположной стороной.
\end{enumerate}
Далее можно убедиться, что ортоцентр, центр тяжести и центр описанной окружности действительно лежат на прямой Эйлера, применив правило \ref{psps} к любой паре из этих точек.

Оказывается описанные правила оригами-построений являются более не менее эфективным инструментом чем циркуль и линейка.
Точнее говоря, любую точку, которую можно построить с помощью циркуля и линейки, можно построить сгибаниями.

Чтобы это доказать, надо предъявить оригами-построения двух типов точек:
\begin{enumerate}
\item Точки пересечения окружности с прямой, если про окружность известно только местоположение центра и одна точка на ней.
\item Точки пересечения двух окружностей, если про окружности известно только местоположение их центров 
и по точке на каждой. 
\end{enumerate}

Действительно, как только мы предъявили нужные оригами-построения, мы можем шаг за шагом дублировать любое построение с помощью циркуля и линейки.

\begin{wrapfigure}{r}{25mm}
\begin{lpic}[t(-8mm),b(-2mm),r(0mm),l(0mm)]{pics/circ-line(1)}
\lbl[t]{12,9;$A$}
\lbl[r]{2,5;$B$}
\lbl[rt]{17,17;$p$}
\end{lpic}
\end{wrapfigure}

Первое можно сделать, применив правило \ref{ptolps}, взяв за $A$ центр окружности, за $B$ точку на окружности, а за $p$ данную прямую а потом применить правило \ref{psls} к $B$ и полученной складке. 
Пересечение второй складки и прямой $p$ есть искомая точка.

Второе сделать сложнее, короткой последовательности сгибаний мне найти не удалось. Но такую последовательность можно получить, показав, что с помощью сгибаний можно построить инверсию точки относительно окружности, если про окружность известно только местоположение центра и одна точка на ней. 
Потом, применив инверсию, которая переводит одну из двух данных окружностей в прямую, свести задачу к предыдущей.

\bigskip
\noi{\bf Упражнение.} \textit{Если вы знаете, что такое инверсия, попробуйте предъявить оригами-построение инверсии точки относительно окружности, если про окружность известно только местоположение центра и одна точка на ней.}
\bigskip

Оказывается, что сгибаниями можно построить точки, которые невозможно построить с помощью циркуля и линейки. 
Вот два примера: эти задачи на построение известны уже более двух тысяч лет, а более сотни лет назад была доказана невозможность решения каждой из этих задач с помощью циркуля и линейки.
Мы предлагаем взять лист бумаги и провести оба построения «руками».

\begin{wrapfigure}{r}{35mm}
\begin{lpic}[t(-8mm),b(-4mm),r(1mm),l(0mm)]{pics/trisect1(1)}
\lbl[rb]{1,3;$A$}
\lbl[b]{31,14;$A'$}
\lbl[lb]{16.5,13.5;$C$}
\lbl[r]{1,22;$B$}
\lbl[l]{21,27;$B'$}
\lbl[t]{27,1;$q$}
\lbl[t]{24,11.5;$q'$}
\lbl[l]{3.5,30;$\ell$}
\lbl[tl]{17,23;$p$}
\end{lpic}
\end{wrapfigure}

\bigskip
\noi{\bf Трисекция угла.} \textit{Разбить данный угол на три равные части.}
\bigskip

\noindent{\it Построение.} 
Пусть угол задан двумя складками $p$ и $q$, обозначим через $A$ вершину угла. 
Сначала проведём подготовительные построения (1) восстановим перпендикуляр $\ell$ к $q$ через $A$ (правило~\ref{psls}) и (2) отметим на $\ell$ произвольно точку $B$ и восстановим срединный перпендикуляр $q'$ к отрезку $AB$ (правило~\ref{ptop}).

Теперь всё готово для главной складки.
Сложим лист так, чтобы $A$ попала на $q'$, а $B$ на $p$ (правило~\ref{2ptol}).
При этом образ $A'$ вершины $A$ ляжет на первую трисектрису нашего угла, 
а точка $C$ на пересечении $q'$ с новой складкой будет лежать на второй. 
То есть лучи $AA'$ и $AC$ будут делить угол на три равные части.

{
\begin{wrapfigure}{r}{36mm}
\begin{lpic}[t(-5mm),b(-3mm),r(0mm),l(0mm)]{pics/kvadratura1(1)}
\lbl[rb]{30.5,2;$A$}
\lbl[rb]{30.5,12;$X$}
\lbl[b]{10,23;$X'$}
\lbl[l]{3,14.5;$A'$}
\lbl[t]{17,20;$q$}
\lbl[r]{1,26;$p$}
\end{lpic}
\end{wrapfigure}

\bigskip
\noi{\bf Удвоение куба.} \textit{Построить два отрезка с отношением длин $\sqrt[3]{2}$.}
\bigskip

\noindent{\it Построение.} Сначала проведём подготовительное построение. 
Построим квадрат, разделённый на три равные части складками параллельными одной стороне. 
Предлагаем это сделать читателю.
Введём обозначения как на чертеже.

}

Теперь сложим лист так, чтобы точка $A$ легла на сторону $p$, a $X$ легла на разделяющую складку $q$. 
При этом $A'$, образ $A$, разделит сторону квадрата в отношении $1:\sqrt[3]2$.

\bigskip
\noi{\bf Упражнение.} \textit{Докажите, что в результате описанных двух построений мы действительно получим то, что хотели.}
\bigskip

На этом список «невозможных построений» не кончается.
Например, при помощи складываний можно также построить правильный семиугольник, который невозможно построить при помощи циркуля и линейки.

В чём же причина, какое из описанных правил добавляет новые возможности? 
Чтобы это понять, можно попытаться построить складку в каждом правиле \ref{psps}---\ref{ptolls}, стр.~\pageref{psps} с помощью циркуля и линейки. 
Довольно легко построить все прямые складок в правилах \ref{psps}---\ref{ptolps} и \ref{ptolls}. 
Стало быть, дополнительные возможности скрыты в правиле номер~\ref{2ptol}. 
Не удивительно что основной шаг в построениях трисекции угла и удвоении куба был сделан применением именно этого правила.
 
Причина, оказывается, в следующем: чтобы найти прямую сгиба в правиле~\ref{2ptol}, требуется решить уравнение третьей степени,
тогда как в каждом из построений с помощью циркуля и линейки решаются уравнения только 1 и 2-ой степеней.

\medskip
\noindent\textbf{Комментарии.}
Заметка написана для моего сына, Никодима.

Правила оригами построений приведённые выше 
обычно называются \emph{правилами Фудзиты}.
Эти семь правил впервые описаны Жаком Жюстином в \cite{justin1989resolution}
чуть позже шесть из них были переоткрыли 
Бенедетто Скимеми и Фумиаки Фудзита \cite{huzita1989algebra},
седьмое правило добавил Хатори Косиро в \cite{hatori2006k}.

Ещё раньше, в \cite{beloch1936sul}, оригами-построения рассматривались Маргеритой Пьяцолла Белок.
Она доказала, что с помощью оригами-построений 
можно решить уравнения третьей и четвёртой степени.
Из этого в частности следует возможность решения задач о трисекции угла и удвоения куба.
О статье Белок я узнал 
из статьи Томаса Хулла \cite{hull},
которую очень рекомендую.

Приведённое нами решение задачи о трисекции угла было предложено  Хисаси Абэ, см. \cite{husimi1980trisection}.
Решение задачи об удвоении куба предложил Петер Мессер в \cite{messer1986problem}.
 
Интересные математические аспекты пространственного оригами обсуждаются в книге Табачнкикова и Фукса \cite{tabach2011}.
В \cite{hull2007constructing}, 
Томас Хулл приводит построения числа $\pi$ с помощью пространственных складок.

Другая известная задача о том что прямоугольник можно сложить на плоскости 
в виде фигуре с б\'{о}льшим периметром 
обсуждается в статье Тарасова \cite{tarasov2004}.
Эта и близкие задачи также обсуждаются в наших лекциях \cite{petrunin2015}.

\end{document}